%% file: the_uniqueness_of_lyapunov_rank_among_symmetric_cones.tex
\documentclass{article}

\usepackage{amsmath, mathtools}

\usepackage{float}

\usepackage[hyphens]{url}

\usepackage{hyperref}
\hypersetup{
  colorlinks=true,
  linkcolor=blue,
  citecolor=blue
}

\usepackage{booktabs}
\input{mjo-algebra}
\input{mjo-arrow}

\input{mjo-common}
\input{mjo-cone}
\input{mjo-font}
\input{mjo-hurwitz}
\input{mjo-proof_by_cases}
\input{mjo-set}
\input{mjo-theorem}

\usepackage[capitalize,nameinlink,noabbrev,nosort]{cleveref}

\crefname{inequality}{Inequality}{Inequalities}
\creflabelformat{inequality}{#2\textup{(#1)}#3}
\crefname{implication}{Implication}{Implications}
\creflabelformat{implication}{#2\textup{(#1)}#3}
\crefname{condition}{Condition}{Conditions}
\creflabelformat{condition}{#2\textup{(#1)}#3}
\crefname{form}{Form}{Forms}
\creflabelformat{form}{#2\textup{(#1)}#3}

\theoremstyle{plain}

\newcommand*{\Roneplus}[0]{ \mathbb{R}^{1}_{+} }

\newcommand*{\signature}[1]{ \sigma\of{#1} }

\newcommand*{\smallestidx}[0]{\iota}

\begin{document}
  \title{
    The uniqueness of Lyapunov rank among symmetric cones
  }
  \author{Michael Orlitzky \and Giovanni Barbarino}

  \maketitle

  \begin{abstract}
    The Lyapunov rank of a cone is the dimension of the Lie algebra of
    its automorphism group. It is invariant under linear isomorphism
    and in general not unique---two or more non-isomorphic cones can
    share the same Lyapunov rank. It is therefore not possible in
    general to identify cones using Lyapunov rank. But suppose we look
    only among \emph{symmetric} cones. Are there any that can be
    uniquely identified (up to isomorphism) by their Lyapunov ranks?
    We provide a complete answer for irreducible cones and make some
    progress in the general case.
  \end{abstract}



  \begin{section}{Introduction}
    The motivation for Lyapunov rank comes from linear programming.
    The familiar complementary slackness condition $\ip{x}{s} = 0$ for
    $x,s \in \Rnplus$ can be separated into $n$ equations
    $\setc{x_{i}s_{i} = 0}{i = 1,2,\ldots,n}$ using the nonnegativity
    of the components, and these $n$ equations should be easier to
    solve simultaneously than the one equation we started
    with. Complementary slackness is a necessary condition for
    optimality, so being able to solve that system can lead us to
    candidate solutions. Rudolf, Noyan, Papp, and
    Alizadeh~\cite{rudolf-noyan-papp-alizadeh} realized that this
    scenario generalizes to a cone program over a proper cone $K$ and
    its dual $\dual{K}$. In a cone program, complementary slackness is
    still a necessary condition for optimality, and the the condition
    $\ip{x}{s} = 0$ for $x \in K$ and $s \in \dual{K}$ can always be
    separated into a system of equations via ``bilinear
    complementarity relations.'' The number of equations, denoted
    $\lyapunovrank{K}$, depends only on the cone $K$ and is called the
    bilinearity rank (by them) or the Lyapunov rank (by us) of $K$. If
    you plan to search for candidate solutions by decomposing the
    complementary slackness condition, Lyapynov rank is a measure of
    how cooperative the cone will be.

    If the cone $K$ lives in an $n$-dimensional space, then the set of
    all candidate pairs $\setc{\pair{x}{s}}{x \in K,\, s \in
      \dual{K},\, \ip{x}{s} = 0}$ forms an $n$-dimensional
    manifold~\cite{rudolf-noyan-papp-alizadeh}. We might therefore
    hope to obtain $n$ or more equations ($\lyapunovrank{K} \ge n$) to
    be solved simultaneously. Rudolf, Noyan, Papp, and Alizadeh call
    such a cone \emph{bilinear}, and show that several important cones
    are not bilinear. A few years later, Gowda and Tao were able to
    connect these ideas to Lyapunov-like
    operators~\cite{gowda-tao_-_rank}. The Lyapunov-like operators on
    $K$ form the Lie algebra of its automorphism group, a vector space
    whose dimension is $\lyapunovrank{K}$. Linearly independent
    Lyapunov-like operators correspond to ``extra'' equations, and
    vice-versa. For this reason the name \emph{Lyapunov rank} was
    chosen.

    Around the same time, Seeger and Sossa coined the term
    \emph{Loewnerian cone} for a cone that is isomorphic to the real
    symmetric PSD cone, $\Snplus$. As demonstrated by Seeger and
    Sossa, Loewnerian cones have important applications to
    complementarity
    problems~\cite{seeger-sossa_-_complementarity_problems}. But how
    do we know if a cone is Loewnerian? Sossa, Gowda, Tao, and one of
    the authors would eventually meet. The problem naturally arose:
    can Lyapunov rank be used to show that a given cone is Loewnerian?
    If the Lyapunov rank of $K$ differs from that of $\Snplus$, then
    $K$ and $\Snplus$ are not isomorphic. The question therefore
    amounts to whether or not there exist any non-Loewnerian $K$
    having $\lyapunovrank{K} = \lyapunovrank{\Snplus}$.

    We are able to answer the question for $\Snplus$ with relative
    ease, but the problem, inspires us to investigate further. First
    we ask whether or not there are any irreducible symmetric
    (self-dual and homogeneous) cones whose Lyapunov ranks are unique
    among other, non-isomorphic symmetric cones. This question is a
    bit harder to resolve, and the answer is unexpected. Afterwards,
    we move on to symmetric cones that are not necessarily
    irreducible. This is related to the isomorphism problem for
    symmetric cones and can have real applications: there do exist
    generic symmetric-cone optimization algorithms, but practically
    speaking, it is much better to know if your symmetric cone is
    really collection of PSD and/or second-order cones in disguise.
  \end{section}

  \begin{section}{Background}
    The main objects of interest are symmetric cones in finite
    dimensions. All such cones arise as the cone of squares in some
    Euclidean Jordan algebra~\cite{faraut-koranyi}. Symmetric cones
    are self-dual, and in particular \emph{proper}: closed, convex,
    pointed, and solid. A convex cone is \emph{reducible} if it is the
    direct sum of two nontrivial convex cones, and \emph{irreducible}
    if not. Every proper cone can be expressed as a direct sum of
    nontrivial irreducible cones, and this sum is unique up to the
    order of its factors (Hauser and Güler~\cite{hauser-guler},
    Theorem~4.3). This general result holds even in the absence of an
    inner product. When we have an inner product, and when our cone is
    symmetric, we can say much more: the decomposition will be
    orthogonal, the factors will be symmetric, and all potential
    factors have been classified up to isomorphism.

    Chapter~V of Faraut and Korányi classifies the simple Euclidean
    Jordan algebras up to isomorphism~\cite{faraut-koranyi}. The
    classification of irreducible symmetric cones then follows from
    the correspondence between symmetric cones and Euclidean Jordan
    algebras. Starting from a symmetric cone, one obtains a Euclidean
    Jordan algebra (Faraut and Korányi, Chapter~III) in which the
    given cone is the cone of squares. That algebra is the orthogonal
    direct sum of simple algebras by Proposition~III.4.4. In each of
    those simple factors, the cone of squares is a symmetric cone, and
    it must be irreducible to avoid contradictions. The resulting
    decomposition of a symmetric cone is Faraut and Korányi's
    Proposition~III.4.5, which we now combine with the classification
    of simple algebras. We use the term \emph{Jordan isomorphism} for
    an invertible Jordan-algebra homomorphism, and
    \emph{Jordan-isomorphic} when two sets are identified by one

    \begin{theorem}\label{thm:direct sum decomposition}
      %
      Every symmetric cone is the unique (up to order) orthogonal
      direct sum of nontrivial irreducible symmetric cones, and every
      nontrivial irreducible symmetric cone is Jordan-isomorphic to a
      member of one of the five distinct families,
      \begin{enumerate}
        \begin{item}
          The $n$-dimensional Lorentz cones, $\Lnplus$, for $n \notin
          \set{0,2}$,
        \end{item}
        \begin{item}
          The $n \times n$ real symmetric positive semidefinite cones,
          $\Hnplus\of{\Rn[1]}$, for $n \ge 3$,
        \end{item}
        \begin{item}
          The $n \times n$ complex Hermitian positive semidefinite
          cones $\Hnplus\of{\Cn[1]}$, for $n \ge 3$,
        \end{item}
        \begin{item}
          The $n \times n$ quaternion Hermitian positive semidefinite
          cones\footnote{These cones of squares were identified as
          positive-semidefinite cones only
          recently~\cite{orlitzky_-_jordan_automorphisms}.},
          $\Hnplus\of{\quaternions}$, for $n \ge 3$,
        \end{item}
        \begin{item}
          The cone of squares in the Albert algebra,
          $\Hnplus[3]\of{\octonions}$.
        \end{item}
      \end{enumerate}
      %
      %
      As a result, every symmetric cone is Jordan-isomorphic to a
      unique (up to order) orthogonal direct sum whose factors come
      from the list above.
    \end{theorem}

    The trivial cone in the trivial space, obtained as an empty sum,
    is also symmetric but is spiritually a member of the family
    $\Lnplus$ with $n=0$. If the ``up to order'' is bothersome, one
    can imagine a total order on the set of irreducible factors listed
    in \cref{thm:direct sum decomposition}. We then obtain true
    uniqueness: every symmetric cone is Jordan-isomorphic to exactly
    one orthogonal direct sum whose factors come from the list and are
    sorted.

    The restrictions on $n$ in \cref{thm:direct sum decomposition}
    ensure that the families are distinct up to Jordan isomorphism.
    The cones $\Hnplus[2]\of{\Rn[1]}$, $\Hnplus[2]\of{\Cn[1]}$,
    $\Hnplus[2]\of{\quaternions}$, and $\Hnplus[2]\of{\octonions}$ are
    symmetric, but they are Jordan-isomorphic to Lorentz cones of
    appropriate sizes (Faraut and Korányi~\cite{faraut-koranyi},
    Corollary~IV.1.5). In fewer than three dimensions, all symmetric
    cones are polyhedral and thus linearly isomorphic. Less obviously,
    they are Jordan-isomorphic: using Proposition~IV.3.2 of Faraut and
    Korányi~\cite{faraut-koranyi}, we can choose bases for both spaces
    consisting of primitive idempotents that belong to extreme rays of
    the symmetric cone. The map sending one basis to the other must
    then be a Jordan isomorphism, and it maps one cone onto the
    other. This is why $\Lnplus[2]$ does not appear in
    \cref{thm:direct sum decomposition}: it is isomorphic to
    $\directsum{\Lnplus[1]}{\Lnplus[1]}$.

    The Lyapunov rank of a proper cone---and in particular of a
    symmetric cone---is the dimension of the Lie algebra of its linear
    automorphism group~\cite{gowda-tao_-_rank}. For our purposes, it
    is simply a number associated with the cone that happens to be
    additive on a direct sum and invariant under linear
    isomorphisms. The following combines Lemma~1 and Proposition~9 of
    Rudolf, Noyan, Papp, and
    Alizadeh~\cite{rudolf-noyan-papp-alizadeh}.

    \begin{proposition}\label{prop:lyapunov rank is additive on direct sums}
      Suppose $K$ and $J$ are symmetric cones, and that $A$ is linear
      and invertible. Then $\lyapunovrank{\directsum{K}{J}} =
      \lyapunovrank{K} + \lyapunovrank{J}$, and $\lyapunovrank{A\of{K}}
        = \lyapunovrank{K}$.
    \end{proposition}

    Taken together, we call the combination of dimension and Lyapunov
    rank the \emph{signature} of a cone.

    \begin{definition}[signature, simulacra]
      Let $K$ be a closed convex cone. The \emph{signature} of $K$ is
      the pair $\signature{K} \coloneqq
      \pair{\dim\of{K}}{\lyapunovrank{K}}$. If $J$ is another closed
      convex cone, we say that $J$ is a \emph{simulacrum} of $K$ if
      $\signature{K} = \signature{J}$ and if $K$ and $J$ are
      non-isomorphic. The plural of simulacrum is \emph{simulacra}. As
      shorthand we write $K \sim J$ to indicate that $K$ and $J$ are
      simulacra of one another. If they are not, then instead we write
      $K \nsim J$.
    \end{definition}

    Gowda and Tao~\cite{gowda-tao_-_rank} computed the Lyapunov ranks
    of the five families in \cref{thm:direct sum decomposition}. Using
    conjugate-symmetry and the fact that
    $\Cn[1],\quaternions,\octonions$ are $2,4,8$-dimensional algebras
    over $\Rn[1]$, it is easy to deduce their dimensions (and thus,
    their signatures) as well. We list these below. Combined with
    \cref{prop:lyapunov rank is additive on direct sums}, and waving
    away the practical difficulty of computing the direct sum
    decomposition in the first place, this gives the signature of any
    symmetric cone.

    \begin{table}[H]
      \begin{center}
        \begin{tabular}{c cc}
          $K$ & $\dim\of{K}$ & $\lyapunovrank{K}$\\
          \midrule
          $\Lnplus$                   & $n$
                                      & $\frac{n^{2} - n + 2}{2}$
          \\
          $\Hnplus\of{\Rn[1]}$        & $\frac{n^{2}+n}2$
                                      & $n^{2}$
          \\
          $\Hnplus\of{\Cn[1]}$        & $n^{2}$
                                      & $2n^{2} - 1$
          \\
          $\Hnplus\of{\quaternions}$  & $2n^{2} - n$
                                      & $4n^{2}$
          \\
          $\Hnplus[3]\of{\octonions}$ & $27$
                                      & $79$
        \end{tabular}
        \caption{Signatures of building blocks, $n \ge 2$}\label{tab:signatures of building blocks}
      \end{center}
    \end{table}

    To this list we add one more, the nonnegative orthant in $n$
    dimensions:

    \begin{table}[H]
      \begin{center}
        \begin{tabular}{c cc}
          $K$ & $\dim\of{K}$ & $\lyapunovrank{K}$\\
          \midrule
          $\Rnplus$ & $n$ & $n$
        \end{tabular}
        \caption{Signature of $\Rnplus$ for $n \ge 0$}\label{tab:rnplus signature}
      \end{center}
    \end{table}

    The cone $\Rnplus \coloneqq \cartprodmany{i=1}{n}{\Lnplus[1]}$ is
    symmetric but not irreducible unless $n \in \set{0,1}$. We include
    it because it is easier, for example, to write $\Rnplus[3]$ than
    it is to write $\Lnplus[1] \oplus \Lnplus[1] \oplus
    \Lnplus[1]$. (The decomposition in \cref{thm:direct sum
      decomposition} is orthogonal, so there is no issue switching
    between Cartesian products and direct sums when we are working up
    to isomorphism.) Keeping in mind that all symmetric cones are
    isomorphic in fewer than three dimensions, we will typically write
    $\Roneplus$ or $\Rnplus[2]$ to indicate a one- or two-dimensional
    symmetric cone. From the two tables, it should be apparent that
    the Lyapunov rank of $\Rnplus$ is strictly less than that of any
    other (non-isomorphic) $n$-dimensional symmetric cone.

    We have alluded to the fact that there are two competing notions
    of isomorphism for a symmetric cone, Jordan and linear. Two convex
    cones $K \subseteq V$ and $J \subseteq W$ are linearly isomorphic
    if there exists an invertible linear $A : V \to W$ such that
    $A\of{K} = J$. This is the usual definition of isomorphism for
    convex cones. In a Euclidean Jordan algebra (that is, for
    symmetric cones) there is the stronger notion of a Jordan algebra
    isomorphism. A Jordan algebra isomorphism is linear and
    invertible, but must also preserve the Jordan algebra
    multiplication. If two Euclidean Jordan algebras are
    Jordan-isomorphic, their cones of squares are linearly isomorphic,
    because the Jordan isomorphism between them preserves squaring and
    is itself a linear isomorphism. We will quickly argue that, for
    symmetric cones, the converse is true as well. This allows us to
    work ``up to isomorphism'' and write things like $K \cong J$ when
    referring to \cref{thm:direct sum decomposition} without having to
    redefine ``isomorphism'' for convex cones.

    \begin{lemma}\label{lem:signatures determine isomorphism for irreducible}
      If $K$ and $J$ are irreducible symmetric cones, then
      $\signature{K} = \signature{J}$ if and only if $K$ and $J$ are
      Jordan-isomorphic.
      \begin{proof}
        Jordan isomorphisms are linear isomorphisms, so one
        implication follows easily from \cref{prop:lyapunov rank is
          additive on direct sums}. For the other, the full list of
        distinct cones (up to Jordan isomorphism) is provided by
        \cref{thm:direct sum decomposition}. According to
        \cref{prop:lyapunov rank is additive on direct sums}, it
        suffices to consider the signatures of only these
        representatives.  With that in mind, \cref{tab:signatures of
          building blocks} can be used to verify that no two distinct
        irreducible factors can have the same signature.
      \end{proof}
    \end{lemma}

    \begin{proposition}
      If $K$ is a symmetric cone in $V$, if $A : V \to W$ is linear
      and invertible, and if $A\of{K}$ is symmetric, then $K$ and
      $A\of{K}$ are Jordan-isomorphic.

      \begin{proof}
        $K$ is symmetric and is therefore the unique (up to order)
        orthogonal direct sum of nontrivial irreducible factors, as in
        \cref{thm:direct sum decomposition},
        \begin{equation*}
          K = K_{1} \oplus K_{2} \oplus \cdots \oplus K_{N}.
        \end{equation*}
        And $A$ is invertible, so it preserves irreducibility: if
        $A\of{K_{i}}$ were expressible as $\directsum{J_{1}}{J_{2}}$,
        we could invert to obtain $K_{i} =
        \directsum{\inverse{A}\of{J_{1}}}{\inverse{A}\of{J_{2}}}$,
        contradicting the irreducibility of $K_{i}$. It follows that
        \begin{equation*}
          A\of{K}
          =
          A\of{K_{1}} \oplus A\of{K_{2}} \oplus \cdots \oplus A\of{K_{N}}
        \end{equation*}
        is a direct sum of nontrivial irreducible proper (but not
        necessarily symmetric) cones. Now, having assumed that
        $A\of{K}$ is symmetric, we may also decompose it into
        nontrivial irreducible symmetric factors using
        \cref{thm:direct sum decomposition}:
        \begin{equation*}
          A\of{K} = J_{1} \oplus J_{2} \oplus \cdots \oplus J_{M}.
        \end{equation*}
        At the beginning of this section we mentioned Theorem~4.3 of
        Hauser and Güler~\cite{hauser-guler}, a weaker version of
        \cref{thm:direct sum decomposition} that can be applied to any
        proper cone. We now cite the uniqueness of this decomposition,
        applied to the proper cone $A\of{K}$, to conclude that our two
        decompositions of $A\of{K}$ have the same factors. As a
        result, each $A\of{K_{i}}$ is some $J_{\ell}$, and is in fact
        symmetric.

        It remains only to show that $A$ does not change the
        Jordan-isomorphism class of any $K_{i}$. Obviously $A$ will
        not change the dimension of $K_{i}$, and \cref{prop:lyapunov
          rank is additive on direct sums} states that $A$ will not
        change its Lyapunov rank either. The question is therefore
        reduced to the following: are there two different (not
        Jordan-isomorphic) irreducible symmetric cones having the same
        signature? This question was answered negatively in
        \cref{lem:signatures determine isomorphism for irreducible}.
      \end{proof}
    \end{proposition}
  \end{section}

  \begin{section}{Irreducible symmetric cones}\label{sec:irreducible cones}
    As a starting point, we show that $\Lnplus$ has no symmetric
    simulacra, regardless of $n$. This is not too surprising; it
    follows from the fact that the Lyapunov rank of $\Lnplus$ is
    strictly maximal in $n$ dimensions. The latter is known to most
    people working with Lyapunov rank, and the details can be worked
    out using \cref{tab:signatures of building blocks}.

    \begin{proposition}\label{prop:rank of lorentz is maximal}
      Among non-isomorphic $n$-dimensional symmetric cones, the
      Lyapunov rank of $\Lnplus$ is strictly maximal.
      \begin{proof}
        First we recall that in dimensions $n \le 2$, all symmetric
        cones are isomorphic. The statement is therefore trivially
        true for $n \le 2$. We next show that splitting $\Lnplus$ into
        two smaller Lorentz cones will necessarily reduce the Lyapunov
        rank. Taking any $k \in \set{1,2,\ldots,n-1}$, and assuming
        now that $n \ge 3$, we can simply compute:
        \begin{equation*}
          \lyapunovrank{\Lnplus}
          -
          \lyapunovrank{\directsum{\Lnplus[k]}{\Lnplus[n-k]} }
          =
          \qty{n-k}{k} - 1.
        \end{equation*}
        For $n \ge 3$, this difference is positive, showing that the
        direct sum has a strictly inferior Lyapunov rank.

        Next, we show that members of the other four irreducible
        families have Lyapunov ranks smaller than a Lorentz cone of
        the same dimension. There's only one cone of octonion matrices
        to worry about, and it's easy to see that
        $\lyapunovrank{\Lnplus[27]} >
        \lyapunovrank{\Hnplus[3]\of{\octonions}}$. For the other three
        families, we compute:
        \begin{alignat*}{3}
          \lyapunovrank{\Lnplus[\qty{m^{2}+m}/2]}
          &-&
          \lyapunovrank{\Hnplus[m]\of{\Rn[1]}}
          &=
          \frac{1}{8}m^{4} + \frac{1}{4}m^{3} - \frac{9}{8}m^{2} - \frac{1}{4}m + 1,
          \\
          \lyapunovrank{\Lnplus[m^{2}]}
          &-&
          \lyapunovrank{\Hnplus[m]\of{\Cn[1]}}
          &=
          \frac{1}{2}m^{4} - \frac{5}{2}m^{2} + 2,
          \\
          \lyapunovrank{\Lnplus[2m^{2} - m]}
          &-&
          \lyapunovrank{\Hnplus[m]\of{\quaternions}}
          &=
          2m^{4} - 2m^{3} - \frac{9}{2}m^{2} + \frac{1}{2}m + 1.
        \end{alignat*}
        At $m = 3$, the higher-order (positive coefficient) terms
        already dominate the lower-order (negative coefficient)
        ones. These expressions therefore remain positive as $m$
        grows. The result now follows: in dimension $n \ge 3$,
        splitting the Lorentz cone into two smaller Lorentz cones will
        decrease the Lyapunov rank, and using anything other than
        Lorentz cones will only decrease it further.
      \end{proof}
    \end{proposition}

    \begin{corollary}\label{cor:lorentz cone has no simulacra}
      $\Lnplus$ has no symmetric simulacra for any $n \ge 0$.
    \end{corollary}

    Next we resolve the problem that prompted this inquiry by showing
    that Lyapunov rank cannot be used to identify Loewnerian cones.

    \begin{proposition}\label{prop:real psd simulacra}
      $\Hnplus\of{\Rn[1]}$ has symmetric simulacra for all $n \ge 3$.
      \begin{proof}
        Compare the signature of $\Lnplus[n+1] \oplus
        \Rnplus[\qty{n^{2} - n - 2}/2]$.
      \end{proof}
    \end{proposition}

    Interestingly, the formula in \cref{prop:real psd simulacra}
    produces a duplicate signature for $n=2$ as well. The resulting
    cone $\Lnplus[3]$ is however isomorphic to $\Hnplus[2]\of{\Rn[1]}$
    and therefore not a simulacrum.

    \begin{proposition}
      $\Hnplus[n]\of{\Cn[1]}$ has symmetric simulacra for $n \ge 4$,
      but not for $n = 3$.

      \begin{proof}
        If $n \ge 4$, then
        \begin{equation*}
          \Hnplus\of{\Cn[1]}
          \sim \Lnplus[n+1]
               \oplus \Lnplus[n+1]
               \oplus \Rnplus[n^{2} - 5n + 6]
               \oplus \Lnplus[4]
               \oplus \underbrace{\Lnplus[3] \oplus \cdots \oplus \Lnplus[3]}_{n-4 \text{ times}}.
        \end{equation*}

        For $n = 3$, we first show that if $\Hnplus[3]\of{\Cn[1]}$ has
        a symmetric simulacrum, then it has one with all Lorentz
        factors---working up to isomorphism, of course. Recall that
        for $m \le 2$, each of the matrix cones
        $\Hnplus[m]\of{\Fn[1]}$ with $\Fn[1] \in \set{\Rn[1], \Cn[1],
          \quaternions, \octonions}$ is in fact a Lorentz cone of an
        appropriate size. We may therefore limit our attention to
        those factors with $m \ge 3$. Having done so,
        $\Hnplus[3]\of{\Rn[1]}$ is the only candidate whose dimension
        allows it to appear (and then only once!) as a factor in a
        simulacrum of $\Hnplus[3]\of{\Cn[1]}$. But if
        $\Hnplus[3]\of{\Cn[1]} \sim
        \directsum{J}{\Hnplus[3]\of{\Rn[1]}}$ for some $J$, then
        \cref{prop:real psd simulacra} tells us that
        $\Hnplus[3]\of{\Cn[1]} \sim J \oplus \Lnplus[4] \oplus
        \Lnplus[2]$ as well. This is a sum of Lorentz cones.

        It now suffices to show that $\Hnplus[3]\of{\Cn[1]}$ does not
        share its signature with a direct sum of Lorentz cones. Since
        $\lyapunovrank{\Lnplus[m]} > 17$ for $m > 6$, we need only
        consider $\Lnplus[m]$ factors with $m \le 6$. And
        $\lyapunovrank{\Lnplus[6]} = 16$, so if $\Lnplus[6]$ were one
        factor, then the other factors would have combined signature
        $\pair{3}{1}$, which is not possible. Thus, $m \le 5$ for all
        $\Lnplus[m]$ factors. If $\Lnplus[5]$ is a factor, then the
        remaining factors have combined signature $\pair{4}{6}$. The
        Lyapunov ranks of $\Rnplus[4]$ and $\Lnplus[3] \oplus
        \Roneplus$ are too small, and that of $\Lnplus[4]$ is too
        large. Thus, a simulacrum cannot have $\Lnplus[5]$ factors.
        If $\Lnplus[4]$ is a factor, then the other factors have
        signature $\pair{5}{10}$. $\Lnplus[5]$ has already been ruled
        out as a factor, and the Lyapunov rank of $\Lnplus[4] \oplus
        \Roneplus$ is too small\@. \cref{prop:rank of lorentz is
          maximal} shows that further splitting of $\Lnplus[4]$ will
        not help, so $\Lnplus[4]$ cannot be a factor. Finally, we note
        that $\lyapunovrank{\Lnplus[3] \oplus \Lnplus[3] \oplus
          \Lnplus[3]}$ is too small.
      \end{proof}
    \end{proposition}

    \begin{proposition}
      $\Hnplus\of{\quaternions}$ has symmetric simulacra for all $n \ge 3$.

      \begin{proof}
        For $n \ge 3$, witness:
        \begin{equation*}
          \Hnplus[n]\of{\quaternions}
          \sim
          \Lnplus[2n+2]
          \oplus
          \Rnplus[2n^{2} - 3n - 2].
          \qedhere
        \end{equation*}
      \end{proof}
    \end{proposition}

    \begin{proposition}\label{prop:octonion psd simulacra}
      $\Hnplus[3]\of{\octonions}$ has symmetric simulacra.

      \begin{proof}
        Compare the signature of $\Lnplus[11] \oplus \Lnplus[5] \oplus
        \Lnplus[3] \oplus \Rnplus[8]$.
      \end{proof}
    \end{proposition}

    Combining all of these, we conclude the section with our first
    major result.

    \begin{theorem}\label{thm:only irreducible factors}
      If $K$ is an irreducible symmetric cone having no symmetric
      simulacra, then either $K \cong \Hnplus[3]\of{\Cn[1]}$ or $K
      \cong \Lnplus$ for some $n \in \Nn[1]$.
    \end{theorem}
  \end{section}

  \begin{section}{Reducible cones}\label{sec:reducible cones}
    The problem becomes much more difficult when we venture beyond
    irreducible cones. We have already met several reducible cones
    whose Lyapunov ranks are not unique---every simulacrum we
    constructed in the previous section has at least one symmetric
    simulacrum, namely the irreducible cone we constructed it to be a
    simulacrum of. So in the general case, the best we can hope for is
    to narrow down the possibilities. In this section we will focus on
    reducing the general problem to smaller subproblems. It is worth
    investigating due to the following (rather obvious) result.

    \begin{lemma}\label{lem:subcone simulacra}
      If any symmetric subcone of a symmetric cone $K$ has symmetric
      simulacra, then so does $K$.
      \begin{proof}
        If $K = \directsum{K_{1}}{K_{2}}$ and if $K_{2} \sim J_{2}$,
        then $K \sim \directsum{K_{1}}{J_{2}}$.
      \end{proof}
    \end{lemma}

    If we can determine which pairs of factors have symmetric
    simulacra, then those pairs cannot appear in any larger symmetric
    cone without inducing symmetric simulacra. This might be called a
    ``bottom up'' approach. Or, suppose we can prove that
    $\directsum{J}{K}$ has symmetric simulacra if and only if $K$
    does. In that case our job has become easier, and if we are lucky,
    the process may be repeated. This would be the ``top down''
    approach.

    \begin{corollary}\label{cor:no multiple complex factors}
      If a symmetric cone contains more than one
      $\Hnplus[3]\of{\Cn[1]}$ factor (up to isomorphism), then it has
      symmetric simulacra.
      \begin{proof}
        If the cone contains two or more $\Hnplus[3]\of{\Cn[1]}$
        factors, then it contains two $\Hnplus[3]\of{\Cn[1]}$ factors,
        and for the purposes of \cref{lem:subcone simulacra} we may use
        \begin{equation*}
          \directsum{\Hnplus[3]\of{\Cn[1]}}{\Hnplus[3]\of{\Cn[1]}}
          \sim
          \Lnplus[7] \oplus \Lnplus[3] \oplus \Rnplus[8].
          \qedhere
        \end{equation*}
      \end{proof}
    \end{corollary}

    \begin{theorem}\label{thm:two sufficient families}
      Up to isomorphism, every symmetric cone $K$ is either of, or has
      simulacra of the form
      \begin{equation*}
        \alpha \Hnplus[3]\of{\Cn[1]}
        \oplus
        \Lnplus[n_{1}]
        \oplus
        \Lnplus[n_{2}]
        \oplus
        \cdots
        \oplus
        \Lnplus[n_{k}]
      \end{equation*}
      where $\alpha \in \set{0,1}$ and $k,n_{1},n_{2},\ldots,n_{k} \in
      \Nn[1]$. If $K$ has no simulacra, then $K$ itself is of this
      form.

      \begin{proof}
        \cref{thm:only irreducible factors} shows that the other
        irreducible families have simulacra involving
        $\Hnplus[3]\of{\Cn[1]}$ and $\Lnplus[n_{i}]$, and \cref{cor:no
          multiple complex factors} shows that we need (or can have)
        at most one $\Hnplus[3]\of{\Cn[1]}$ factor.
      \end{proof}
    \end{theorem}

    This result can be interpreted as saying that every symmetric cone
    shares its signature with a direct sum of Lorentz cones and/or
    $\Hnplus[3]\of{\Cn[1]}$. We caution however that simulacra of this
    form may not be unique. Consider the simulacra,
    \begin{equation*}
      \Lnplus[4] \oplus \Lnplus[3] \oplus \Lnplus[3] \oplus \Lnplus[3]
      \sim
      \Lnplus[5] \oplus \Rnplus[8].
    \end{equation*}
    This shows that a cone can be isomorphic to, \emph{and} have
    simulacra of the form in \cref{thm:two sufficient
      families}. Moreover, when applied to \cref{prop:octonion psd
      simulacra}, it shows that one cone can have two different
    simulacra of the desired form.

    \cref{thm:two sufficient families} suggests which irreducible
    factors may be interesting to pair up and analyze. Eventually we
    will look at both $\directsum{\Hnplus[3]\of{\Cn[1]}}{\Lnplus}$ and
    $\directsum{\Lnplus[m]}{\Lnplus[n]}$, but to lay the groundwork,
    we first consider the general case of $\directsum{K}{\Lnplus}$. If
    the problem is to find simulacra of $\directsum{K}{\Lnplus}$, we
    will impose three conditions on $n$,
    \begin{enumerate}
      \begin{item}\label[condition]{itm:condition1}
        $n \ge 2 + \lyapunovrank{K} - \dim\of{K}$,
      \end{item}
      \begin{item}\label[condition]{itm:condition2}
        $n \ge 2 + \lyapunovrank{\Lnplus[1 + \dim\of{K}]} - \lyapunovrank{K}$,
      \end{item}
      \begin{item}\label[condition]{itm:condition3}
        $n \ge 15$.
      \end{item}
    \end{enumerate}
    The $n \ge 15$ bound is not very elegant, but it will simplify our
    proofs. Either way we are going to wind up checking cases on a
    computer; but this way, the computer does a bit more work and we
    do a little less.

    \begin{lemma}\label{lem:condition4}
      If $K$ is a symmetric cone and if $n \in \Nn[1]$ satisfies
      \cref{itm:condition1,itm:condition2,itm:condition3}, then $n >
      2\dim\of{K}$.

      \begin{proof}
        From the assumption that $n \ge 15$, the conclusion is
        obviously true whenever $\dim\of{K} \le 7$. Averaging
        \cref{itm:condition1,itm:condition2}, we obtain a quadratic
        function of $\dim\of{K}$ that serves as another lower bound on
        $n$. One can set this new bound to be greater than
        $2\dim\of{K}$ and solve the quadratic to see that
        $\dim\of{K} \ge 8$ implies $n > 2\dim\of{K}$.
      \end{proof}
    \end{lemma}

    We next use \cref{itm:condition1} to rule out the possibility that
    $\directsum{K}{\Lnplus}$ shares its signature with some other cone
    having a larger Lorentz cone factor than $\Lnplus$.

    \begin{lemma}\label{lem:no bigger factors}
      If $K$ is a symmetric cone and if $n \in \Nn[1]$ satisfies
      \cref{itm:condition1}, then
      $\signature{\directsum{J}{\Lnplus[n+k]}} \ne
      \signature{\directsum{K}{\Lnplus}}$ for any symmetric cone $J$
      and any $k \ge 1$.
      \begin{proof}
        With $k \ge 1$, the smallest possible Lyapunov rank achievable
        by $\directsum{J}{\Lnplus[n+k]}$ in dimension $\dim\of{K}+n$
        is at $k = 1$ with $J = \Rnplus[\dim\of{K} - k]$. But if $n$
        satisfies \cref{itm:condition1}, then
        \begin{equation*}
          \lyapunovrank{\Rnplus[\dim\of{K} - 1] \oplus \Lnplus[n+1]}
          -
          \lyapunovrank{K \oplus \Lnplus}
          =
          \dim\of{K} - 1 - \lyapunovrank{K} + n
          \ge
          1.
        \end{equation*}
        Thus, the smallest possible Lyapunov rank is too large to work.
      \end{proof}
    \end{lemma}

    The other conditions can now be used to rule out sums of Lorentz
    cones wherein each factor is smaller than $\Lnplus$.

    \begin{lemma}\label{lem:not all smaller factors}
      Suppose $K$ is a symmetric cone, that $n \in \Nn[1]$ satisfies
      \cref{itm:condition1,itm:condition2,itm:condition3}, and that $k
      \in \Nn[1]$. Then
      \begin{equation*}
        \lyapunovrank{\directsummany{i=1}{k}{\Lnplus[n_{i}]}}
        <
        \lyapunovrank{\directsum{K}{\Lnplus}}
      \end{equation*}
      for all $n_{1},n_{2},\ldots,n_{k} < n$ such that
      $\dim\of{\directsummany{i=1}{k}{\Lnplus[n_{i}]}} =
      \dim\of{\directsum{K}{\Lnplus}}$.

      \begin{proof}
        Note that $k \ge 2$ if $n_{1} < n$. By rearranging the
        factors, we may assume that $n_{1} \ge n_{2} \ge \cdots \ge
        n_{k}$. To simplify the notation, we let
        \begin{equation*}
          f\of{x} \coloneqq \lyapunovrank{\Lnplus[x]},
        \end{equation*}
        and notice that, by \cref{prop:rank of lorentz is maximal},
        \begin{equation}\label[inequality]{ineq:f is superadditive}
          f\of{x} + f\of{y} \le f\of{x + y}.
        \end{equation}
        It is also easy to see for $n \ge 2$, which is guaranteed by
        \cref{itm:condition3}, that
        \begin{equation}\label{eq:f of n - 1}
          f\of{n - 1} = f\of{n} - \qty{n - 1}.
        \end{equation}
        If we fix $\delta \ge 0$, then $g_{\delta}\of{x} \coloneqq
        f\of{x} - f\of{x - \delta}$ is nondecreasing since
        $g_{\delta}'\of{x} = \delta$. In particular, $g_{\delta}\of{n
          - 1} \ge g_{\delta}\of{1 + \dim\of{K} + \delta}$ whenever $n
        - 1 \ge 1 + \dim\of{K} + \delta$. In other words,
        \begin{equation}\label[implication]{impl:growth inequality for f pt2}
          \begin{gathered}
            n - 1 \ge 1 + \dim\of{K} + \delta\\
            \implies\\
            f\of{n - 1} + f\of{1 + \dim\of{K}}
            \ge
            f\of{n - 1 - \delta} + f\of{1 + \dim\of{K} + \delta}.
          \end{gathered}
        \end{equation}

        Let $S\of{j} \coloneqq \sum_{i=1}^{j}n_{i}$, and define
        $\smallestidx$ to be the smallest index such that
        $S\of{\smallestidx} > \dim\of{K}$. Since $S\of{k} = n +
        \dim\of{K} > \dim\of{K}$, this is well-defined. From
        \cref{lem:condition4} we have that $n > \dim\of{K}$,
        which moreover implies that $\smallestidx < k$. Applying
        \cref{ineq:f is superadditive} repeatedly, we find that
        \begin{equation*}
          \sum_{i=1}^{k}f\of{n_{i}}
          =
          \sum_{i=1}^{\smallestidx}f\of{n_{i}}
          +
          \sum_{i=\smallestidx+1}^{k}f\of{n_{i}}
          \le
          f\of{S\of{\smallestidx}}
            + f\of{n + \dim\of{K}
            - S\of{\smallestidx}}.
        \end{equation*}
        Let $\delta \coloneqq S\of{\smallestidx} - \dim\of{K} -
        1$. Then $\delta \ge 0$ by the definition of $\smallestidx$,
        and $\dim\of{K} + 1 + \delta = S\of{\smallestidx}$. To apply
        \cref{impl:growth inequality for f pt2}, we would like to know
        that $S\of{\smallestidx} \le n - 1$. Using
        \cref{lem:condition4}: if $n_{1} > \dim\of{K}$, then
        $\smallestidx = 1$ and $S\of{\smallestidx} = n_{1} < n$ by
        assumption; if $n_{1} \le \dim\of{K}$, then
        $S\of{\smallestidx} = S\of{\smallestidx - 1} +
        n_{\smallestidx} \le 2\dim\of{K} < n$.

        We may now replace $S\of{\smallestidx}$ by $\dim\of{K} + 1 +
        \delta$ and apply \cref{impl:growth inequality for f
          pt2}:
        \begin{align*}
          \sum_{i=1}^{k}f\of{n_{i}}
          &\le
          f\of{1 + \dim\of{K} + \delta} + f\of{n - 1 - \delta}\\
          &\le
          f\of{n-1} + f\of{1 + \dim\of{K}}.
        \end{align*}
        Finally, using \cref{eq:f of n - 1}, we arrive at
        \begin{equation*}
          \sum_{i=1}^{k}f\of{n_{i}}
          \le
          f\of{n} - \qty{n-1}+ f\of{1 + \dim\of{K}},
        \end{equation*}
        which, by \cref{itm:condition2}, is strictly less than
        $f\of{n} + \lyapunovrank{K}$.
      \end{proof}
    \end{lemma}

    \begin{theorem}\label{thm:subproblem reduction}
      Suppose that $K,J$ are symmetric cones and that $n \in \Nn[1]$
      satisfies \cref{itm:condition1,itm:condition2,itm:condition3}.
      Then $J \sim \directsum{K}{\Lnplus}$ if and only if there exists
      a symmetric $J'$ such that $J \cong \directsum{J'}{\Lnplus}$ and
      $J' \sim K$.

      \begin{proof}
        The ``if'' direction and the case $K = \set{0}$ are obvious
        and don't require any conditions on $n$. In the other
        direction, the main difficulty is showing that $J$ must itself
        have an $\Lnplus$ factor. To that end, suppose that $J \sim
        \directsum{K}{\Lnplus}$. We conclude from \cref{lem:no bigger
          factors} that none of the Lorentz factors in $J$ are larger
        than $\Lnplus$, and \cref{lem:not all smaller factors}
        moreover implies that $J$ cannot consist entirely of Lorentz
        cones of dimension $n-1$ or less. It follows that if $J$ has
        only Lorentz factors, then it has at least one $\Lnplus$
        factor.

        If $J$ \emph{does} contain non-Lorentz factors, we can
        ``replace'' each non-Lorentz factor $I$ by a sum of Lorentz
        cones, all of dimension $n-1$ or less, to obtain a new
        symmetric cone $H$. The intention is that this $H$ will have
        the same dimension as $J$ and will satisfy $\lyapunovrank{J}
        \le \lyapunovrank{H} < \lyapunovrank{\directsum{K}{\Lnplus}}$
        per \cref{lem:not all smaller factors} unless $\Lnplus$ is a
        factor of $J$. Refer back to \cref{tab:signatures of building
          blocks}: if $J$ contains an irreducible, non-Lorentz factor
        $I$, then $I$ satisfies $\dim\of{I} \ge \lyapunovrank{I}/3$.
        And either $I \cong \Hnplus[3]\of{\Cn[1]}$, or, using the
        results in \cref{sec:irreducible cones}, $I$ has a simulacrum
        $I'$ all of whose factors are Lorentz cones. Now, consider
        each non-Lorentz factor $I$ of $J$, one at a time. If $I \cong
        \Hnplus[3]\of{\Cn[1]}$, replace it by $\Lnplus[9]$, and note
        that $9 < n$ by \cref{itm:condition3}. If $I \ncong
        \Hnplus[3]\of{\Cn[1]}$, replace it by its all-Lorentz
        simulacrum $I'$. We claim that none of the Lorentz factors in
        $I'$ are of dimension $n$ or greater. If $I'$ contains an
        $\Lnplus[n+k]$ factor where $k \ge 0$, then
        \begin{equation*}
          \dim\of{\directsum{\Lnplus}{K}}
          \ge
          \dim\of{I'}
          \ge
          \frac{1}{3}\lyapunovrank{I'}
          \ge
          \frac{1}{3}\lyapunovrank{\Lnplus[n+k]}
          \ge
          \frac{1}{3}\lyapunovrank{\Lnplus[n]}
          \ge
          2n,
        \end{equation*}
        where the last inequality is due to
        \cref{itm:condition3}. Subtracting $n$ from both sides we
        would obtain $\dim\of{K} \ge n$ in contradiction of
        \cref{lem:condition4}. Thus all Lorentz factors of $I'$ are of
        dimension $n-1$ or less. Proceeding in this manner, we
        construct the cone $H$ which consists of:
        \begin{itemize}
          \begin{item}
            The Lorentz factors in $J$,
          \end{item}
          \begin{item}
            An $\Lnplus[9]$ for every $I \cong \Hnplus[3]\of{\Cn[1]}$
            in $J$,
          \end{item}
          \begin{item}
            A direct sum of Lorentz cones, all of dimension $n-1$ or less,
            for every non-Lorentz $I \ncong \Hnplus[3]\of{\Cn[1]}$ in $J$.
          \end{item}
        \end{itemize}
        And, by construction, $\dim\of{H} = \dim\of{J}$ with
        $\lyapunovrank{H} \ge \lyapunovrank{J}$. Recall that $J$
        contains no Lorentz factors larger than $\Lnplus$ by
        \cref{lem:no bigger factors}. If $J$ also contains no
        $\Lnplus$ factors, then all factors of $H$ are Lorentz cones
        of dimension $n-1$ or less. In that case we would conclude
        from \cref{lem:not all smaller factors} that $\lyapunovrank{J}
        \le \lyapunovrank{H} < \lyapunovrank{\directsum{K}{\Lnplus}}$,
        but of course this is not possible. The only other option is
        for $J$ to have an $\Lnplus$ factor.

        Now that we know it exists, take one of the $\Lnplus$ factors
        from within $J$, and write $J \cong
        \directsum{J'}{\Lnplus}$. It should be clear that if $J \sim
        \directsum{K}{\Lnplus}$, then the signatures of $J'$ and $K$
        agree. Moreover, $J'$ and $K$ cannot be isomorphic; for if
        they were, then $J$ and $\directsum{K}{\Lnplus}$ would be
        isomorphic, and they are not.
      \end{proof}
    \end{theorem}

    \cref{itm:condition3} can be relaxed so long as we keep the other
    two. Adding \cref{itm:condition1,itm:condition2}, we obtain the
    inequality $\dim\of{K}\qty{\dim\of{K} - 1} \le 4n - 10$. For this
    inequality to hold we must have $n \ge 3$, and the potential new
    cases are $n \in \set{3,4,\ldots,14}$ because $n \ge 15$ is
    handled by the theorem. For each $n$, then, we must check all $K$
    such that $\dim\of{K}\qty{\dim\of{K} - 1} \le 4n - 10$. There are
    a finite number of these, and computation shows that
    $\Hnplus[3]\of{\Rn[1]} \sim \directsum{\Lnplus[4]}{\Rnplus[2]}$ is
    the sole counterexample at $n = 4$. As a result, we can get away
    with requiring $n \ne 4$ instead of $n \ge 15$.

    In the proof of this result, we will use integer partitioning to
    enumerate sums of Lorentz cones, a strategy we employ several
    times. There are fast algorithms to compute partitions; for
    example, the accelerated integer partitioning scheme of Kelleher
    and O'Sullivan~\cite{kelleher-osullivan}. Most partitioning
    schemes sort the entries to avoid duplication, so no ambiguity
    arises in the order of our Lorentz factors. The partitions $1 + 1$
    and $2$ however correspond to isomorphic cones---we do not need to
    consider both. Fortunately, as the entries of each partition are
    already sorted, it is easy to detect when a partition contains $2$
    and omit it from the result. For every partition containing $2$
    that we omit, there will be a partition with $1+1$ in its place,
    so we do not lose any cones. This can save a considerable amount
    of time.

    \begin{theorem}[Improved \cref{thm:subproblem reduction}]\label{thm:improved subproblem reduction}
      Suppose that $K,J$ are symmetric cones and that $n \ne 4$
      satisfies \cref{itm:condition1,itm:condition2}. Then $J \sim
      \directsum{K}{\Lnplus}$ if and only if there exists a symmetric
      $J'$ such that $J \cong \directsum{J'}{\Lnplus}$ and $J' \sim
      K$.

      \begin{proof}
        Recall from the proof of \cref{thm:subproblem reduction} that
        the ``if'' direction was obvious and required no conditions on
        $n$. For the other direction, as just discussed, we must check
        the cones $K$ satisfying $\dim\of{K}\qty{\dim\of{K} - 1} \le
        4n - 10$ for each $n \in \set{3,5,6,\ldots,14}$. We begin by
        computing the maximum value of $\dim\of{K}$ corresponding to
        each $n$.
        \begin{table}[H]
          \begin{center}
            \begin{tabular}{l l}
              $n$ & $\dim\of{K}$\\
              \midrule
              $3$ & $1,2$\\
              $5$ & $1,2,3$\\
              $6,7$ & $1,2,3,4$\\
              $8,9$ & $1,2,3,4,5$\\
              $10,11,12$ & $1,2,3,4,5,6$\\
              $13,14$ & $1,2,3,4,5,6,7$\\
            \end{tabular}
          \end{center}
        \end{table}
        The total dimension of $J$ is therefore at most $21$. This
        will allow us to eliminate all but one non-Lorentz cone as
        possible factors of $J$ and $K$. The only non-Lorentz,
        irreducible symmetric cone of dimension $7$ or less is
        $\Hnplus[3]\of{\Rn[1]}$, so immediately we see that no other
        non-Lorentz factors can be present in $K$. If we allow
        dimensions of $21$ or less (for factors of $J$), then a few
        more possibilities arise: $\Hnplus[4]\of{\Rn[1]}$,
        $\Hnplus[5]\of{\Rn[1]}$, $\Hnplus[6]\of{\Rn[1]}$,
        $\Hnplus[3]\of{\Cn[1]}$, $\Hnplus[4]\of{\Cn[1]}$, and
        $\Hnplus[3]\of{\quaternions}$. However, each of these can be
        ruled out.

        Recall that $\lyapunovrank{\directsum{K}{\Lnplus}}$ is bounded
        below by $\lyapunovrank{\Lnplus} + \dim\of{K}$, since $K =
        \Rnplus[\dim\of{K}]$ minimizes the Lyapunov rank. If $J$ has
        an $\Hnplus[3]\of{\Cn[1]}$ factor, then $\lyapunovrank{J} = 17
        + \lyapunovrank{J'} \le 17 + \lyapunovrank{\Lnplus[n +
            \dim\of{K} - 9]}$ by \cref{prop:rank of lorentz is
          maximal}. Subtracting and substituting the pairs of $n$ and
        $\dim\of{K}$ from the table above, we find that
        \begin{equation*}
          \lyapunovrank{\Lnplus} + \dim\of{K}
          -17
          -\lyapunovrank{\Lnplus[n + \dim\of{K} - 9]}
          > 0,
        \end{equation*}
        for all compatible $n$ and $\dim\of{K}$ such that $n +
        \dim\of{K} \ge 9$; that is, all pairs that would make $J$
        large enough to hold an $\Hnplus[3]\of{\Cn[1]}$ factor. An
        analogous argument rules out $\Hnplus[4]\of{\Rn[1]}$,
        $\Hnplus[4]\of{\Cn[1]}$, and
        $\Hnplus[3]\of{\quaternions}$. For the remaining two, we can
        take a shortcut by noticing that
        \begin{align*}
          \Hnplus[5]\of{\Rn[1]}
          &\sim
          \Hnplus[3]\of{\Cn[1]} \oplus \Lnplus[3] \oplus \Lnplus[3],\\
          \Hnplus[6]\of{\Rn[1]}
          &\sim
          \Hnplus[3]\of{\Cn[1]}
            \oplus \Lnplus[4]
            \oplus \Lnplus[4]
            \oplus \Lnplus[3]
            \oplus \Roneplus.
        \end{align*}
        The $\Hnplus[3]\of{\Cn[1]}$ factor in these makes their
        Lyapunov ranks too small, so neither of them can be factors in
        $J$. This leaves only $\Hnplus[3]\of{\Rn[1]}$ to consider. The
        strategy we have used so far cannot rule
        $\Hnplus[3]\of{\Rn[1]}$ out entirely. We may however draw
        several conclusions:
        \begin{itemize}
          \begin{item}
            $\Hnplus[3]\of{\Rn[1]}$ cannot be a factor of $K$, because
            if it were, then $K$ would be isomorphic to either
            $\Hnplus[3]\of{\Rn[1]}$ or
            $\directsum{\Hnplus[3]\of{\Rn[1]}}{\Roneplus}$. In either
            of those cases, \cref{itm:condition2} requires $n \ge 15$.
          \end{item}

          \begin{item}
            $\Hnplus[3]\of{\Rn[1]}$ cannot be a factor of $J$ unless
            $n \ge 10$ and $\dim\of{K} \ge 6$ because, otherwise, the
            Lyapunov rank of $J$ will be too small (the same argument
            we used a moment ago).
          \end{item}

          \begin{item}
            $\directsum{\Hnplus[3]\of{\Rn[1]}}{\Hnplus[3]\of{\Rn[1]}}$
            cannot appear as a factor of $J$ for the same reason.
          \end{item}
        \end{itemize}
        So at most we can have a single $\Hnplus[3]\of{\Rn[1]}$ factor
        in $J$, and only when $n \ge 10$. For lack of a better
        approach, the remaining cases can be checked by brute
        force. If we had only Lorentz factors in both $K$ and $J$, we
        could enumerate the possibilities by partitioning the integers
        $\dim\of{K}$ and $n + \dim\of{K}$ for $K$ and $J$
        respectively. But as there may be an $\Hnplus[3]\of{\Rn[1]}$
        in $J$, the process gets a little more complicated.

        For $n \le 9$ a pure integer-partition comparison can be used.
        For each $n \ge 10$ and each associated $\dim\of{K}$ in the
        table, there are a two cases. For $K$, we have the integer
        partitions of $\dim\of{K}$; but for $J$, we have not only the
        sums of Lorentz cones arising from the integer partitions of
        $n + \dim\of{K}$, but also those cones that arise from a
        direct sum of $\Hnplus[3]\of{\Rn[1]}$ and an integer partition
        of $n + \dim\of{K} - 6$. Essentially we decide whether or not
        to include $\Hnplus[3]\of{\Rn[1]}$ in $J$, and then partition
        the remaining dimensions.

        Once the possibilities are computed for a fixed $n$ and
        $\dim\of{K}$, we compare Lyapunov ranks. Keeping in mind that
        \cref{itm:condition1,itm:condition2} must be satisfied, it
        easy to see that the theorem holds because $J \nsim
        \directsum{K}{\Lnplus}$ for all valid $K,J$.
      \end{proof}
    \end{theorem}

    None of the remaining conditions can be
    weakened\@. \cref{itm:condition1} was chosen to exceed the value
    of an explicit counterexample in \cref{lem:no bigger factors}. If
    $n = 1 + \lyapunovrank{K} - \dim\of{K}$, one fewer than the bound,
    then excepting a few pathological $K$,
    \begin{equation*}
      \directsum{K}{\Lnplus}
      \sim
      \directsum{\Rnplus[\dim\of{K}-1]}{\Lnplus[n+1]}.
    \end{equation*}
    By taking $K \coloneqq \Lnplus[m]$ for $m \ge 5$, we can be sure
    that this remains a counterexample even when the other two
    conditions are satisfied. When $K = \Rnplus[2]$ and $n = 4$, only
    the condition $n \ne 4$ is violated, and we get the simulacra
    $\Hnplus[3]\of{\Rn[1]} \sim \directsum{\Rnplus[2]}{\Lnplus[4]}$.
    Finally, when $K = \Hnplus[3]\of{\Cn[1]}$, \cref{itm:condition2}
    gives the bound $n \ge 31$. But for $n = 30$, the other two
    conditions are satisfied, and
    \begin{equation*}
      \directsum{\Hnplus[3]\of{\Cn[1]}}{\Lnplus}
      \sim
      \directsum{\Lnplus[29]}{\Lnplus[10]}.
    \end{equation*}
    Each condition is therefore essential.

    \begin{corollary}
      If $K$ is a symmetric cone having no symmetric simulacra and if
      $n \ne 4$ satisfies \cref{itm:condition1,itm:condition2}, then
      $\directsum{K}{\Lnplus}$ has no symmetric simulacra. In
      particular, $\directsum{\Hnplus[3]\of{\Cn[1]}}{\Lnplus}$ has no
      symmetric simulacra for $n \ge 31$.
    \end{corollary}

    To catalogue the simulacra of
    $\directsum{\Hnplus[3]\of{\Cn[1]}}{\Lnplus}$, it now suffices to
    check the first thirty values of $n$. And we need only consider
    simulacra the form $\directsummany{i=1}{k}{\Lnplus[n_{i}]}$ with
    $n-1 \ge n_{1} \ge n_{2} \ge \cdots \ge n_{k}$ and $k \ge 2$. For
    if a simulacrum contains $\Hnplus[3]\of{\Cn[1]}$, then, by
    subtracting $\pair{9}{17}$ from both signatures, it also contains
    a simulacrum for $\Lnplus$, which is not possible by
    \cref{cor:lorentz cone has no simulacra}. And if a simulacrum does
    not contain $\Hnplus[3]\of{\Cn[1]}$, then either it consists of
    only Lorentz factors, or the results in \cref{sec:irreducible
      cones} let us find another cone having the same signature and
    only Lorentz factors. Since we need only consider Lorentz factors,
    enumerating the possibilities is easily accomplished by
    partitioning $9+n$. Absent rows indicate the absence of simulacra.
    \begin{table}[H]
      \begin{center}
        \begin{tabular}{c ccc}
          $K$ & $\dim\of{K}$ & $\lyapunovrank{K}$ & $\sim$\\
          \midrule
          $\directsum{\Hnplus[3]\of{\Cn[1]}}{\Lnplus[2]}$ &
            $11$ &
            $19$ &
            $\Lnplus[5] \oplus \Lnplus[3] \oplus \Lnplus[3]$\\
          $\directsum{\Hnplus[3]\of{\Cn[1]}}{\Lnplus[3]}$ &
            $12$ &
            $21$ &
            $\Lnplus[4] \oplus \Lnplus[4] \oplus \Lnplus[4]$\\
          $\directsum{\Hnplus[3]\of{\Cn[1]}}{\Lnplus[4]}$ &
            $13$ &
            $24$ &
            $\Lnplus[6] \oplus \Lnplus[3] \oplus \Rnplus[4]$\\
          $\directsum{\Hnplus[3]\of{\Cn[1]}}{\Lnplus[5]}$ &
            $14$ &
            $28$ &
            $\Lnplus[6] \oplus \Lnplus[4] \oplus \Lnplus[3] \oplus \Roneplus$\\
          $\directsum{\Hnplus[3]\of{\Cn[1]}}{\Lnplus[6]}$ &
            $15$ &
            $33$ &
            $\Lnplus[5] \oplus \Lnplus[5] \oplus \Lnplus[5]$\\
          $\directsum{\Hnplus[3]\of{\Cn[1]}}{\Lnplus[7]}$ &
            $16$ &
            $39$ &
            $\Lnplus[6] \oplus \Lnplus[6] \oplus \Lnplus[4]$\\
          $\directsum{\Hnplus[3]\of{\Cn[1]}}{\Lnplus[8]}$ &
            $17$ &
            $46$ &
            $\Lnplus[9] \oplus \Lnplus[3] \oplus \Rnplus[5]$\\
          $\directsum{\Hnplus[3]\of{\Cn[1]}}{\Lnplus[9]}$ &
            $18$ &
            $54$ &
            $\Lnplus[10] \oplus \Rnplus[8]$\\
          $\directsum{\Hnplus[3]\of{\Cn[1]}}{\Lnplus[10]}$ &
            $19$ &
            $63$ &
          $\Lnplus[9] \oplus \Lnplus[7] \oplus \Lnplus[3]$\\
          $\directsum{\Hnplus[3]\of{\Cn[1]}}{\Lnplus[15]}$ &
            $24$ &
            $123$ &
          $\Lnplus[14] \oplus \Lnplus[8] \oplus \Rnplus[2]$\\
          $\directsum{\Hnplus[3]\of{\Cn[1]}}{\Lnplus[18]}$ &
            $27$ &
            $171$ &
            $\Lnplus[14] \oplus \Lnplus[13]$\\
          $\directsum{\Hnplus[3]\of{\Cn[1]}}{\Lnplus[21]}$ &
            $30$ &
            $228$ &
            $\Lnplus[19] \oplus \Lnplus[11]$\\
          $\directsum{\Hnplus[3]\of{\Cn[1]}}{\Lnplus[22]}$ &
            $31$ &
            $249$ &
            $\Lnplus[21] \oplus \Lnplus[9] \oplus \Roneplus$\\
          $\directsum{\Hnplus[3]\of{\Cn[1]}}{\Lnplus[30]}$ &
            $39$ &
            $453$ &
            $\Lnplus[29] \oplus \Lnplus[10]$\\
          \end{tabular}
        \caption{Simulacra of $\directsum{\Hnplus[3]\of{\Cn[1]}}{\Lnplus}$}
      \end{center}
    \end{table}

    \begin{proposition}
      $\Hnplus[3]\of{\Cn[1]} \oplus \Lnplus$ has symmetric simulacra
      if and only if $n \in \set{2,3,\ldots,10,15,18,21,22,30}$. If
      $K$ is a symmetric cone with an $\Hnplus[3]\of{\Cn[1]}$ factor
      and no symmetric simulacra, then
      \begin{equation*}
        K
        \cong
        \Hnplus[3]\of{\Cn[1]}
        \oplus
        \Lnplus[n_{1}]
        \oplus
        \Lnplus[n_{2}]
        \oplus
        \cdots
        \oplus
        \Lnplus[n_{k}]
      \end{equation*}
      where $k \in \Nn[1]$ and $n_{i} \in \Nn[1]
      \setminus \set{2,3,\ldots,10,15,18,21,22,30}$.
      \begin{proof}
        With \cref{lem:subcone simulacra} in mind, start from
        \cref{thm:two sufficient families} and remove the known simulacra
        from the table.
      \end{proof}
    \end{proposition}

    \begin{proposition}
      If $m \ne 2$ and if $n \ge \qty{m^2 - 3m + 6}/2$, then
      $\directsum{\Lnplus[m]}{\Lnplus}$ has no symmetric
      simulacra. When $m = 2$, the only counterexample that arises is
      $\directsum{\Rnplus[2]}{\Lnplus[4]} \sim \Hnplus[3]\of{\Rn[1]}$.
      \begin{proof}
        The $m = 0$ case is handled by \cref{cor:lorentz cone has no
          simulacra}, so we may assume that $m \ge 1$. The condition
        $n \ge \qty{m^2 - 3m + 6}/2$ is thus \cref{itm:condition1}
        applied to $K \coloneqq \Lnplus[m]$.

        When $m \ge 5$, \cref{itm:condition1} dominates
        \cref{itm:condition2} for $K \coloneqq \Lnplus[m]$ and
        moreover implies that $n \ge 5$, so this is a corollary of
        \cref{thm:improved subproblem reduction} in that case. The
        remaining pairs of $m,n$ that satisfy $m \ge 1$ and $n \ge
        \qty{m^2 - 3m + 6}/2$ but \emph{not} the prerequisites for
        \cref{thm:improved subproblem reduction} are,

        \begin{itemize}
          \begin{item}
            $m = 4$, where \cref{itm:condition2} is violated for $n = 5$.
          \end{item}
          \begin{item}
            $m = 3$, where \cref{itm:condition2} is violated for $n
            \in \set{3,4}$.
          \end{item}
          \begin{item}
            $m = 2$, where $n \ne 4$ is violated for $n = 4$ and
            \cref{itm:condition2} is violated for $n \in \set{2,3}$.
          \end{item}
          \begin{item}
            $m = 1$, where \cref{itm:condition2} is violated for $n =
            2$ and $n \ne 4$ is violated for $n = 4$.
          \end{item}
        \end{itemize}
        All of these cases are of dimension nine or less, where
        $\Hnplus[3]\of{\Rn[1]}$ is the only non-Lorentz cone small
        enough to appear in a simulacrum. A priori, to check them all,
        we would need to check all sums of Lorentz cones with/without
        an $\Hnplus[3]\of{\Rn[1]}$ factor. But there are some
        shortcuts we can take:

        \begin{itemize}
          \begin{item}
            $\directsum{\Lnplus[m]}{\Lnplus} \cong \Rnplus[m+n]$ has
            no simulacra when $m,n \le 2$.
          \end{item}
          \begin{item}
            $\directsum{\Lnplus[3]}{\Lnplus[2]}$ and
            $\directsum{\Lnplus[4]}{\Lnplus[1]}$ have no
            simulacra. Neither $\Lnplus[5]$ nor any $\Rnplus[k]$ has
            simulacra, so they would have to be simulacra of each
            other (they aren't).
          \end{item}
          \begin{item}
            When $m,n \ge 3$, there are only two factors in
            $\directsum{\Lnplus[m]}{\Lnplus}$. As a result it cannot
            be isomorphic to any symmetric cone whose decomposition
            contains three or more factors; if
            $\directsum{\Lnplus[m]}{\Lnplus}$ shares its signature
            with such a cone, they are simulacra. From \cref{prop:real
              psd simulacra}, we know that $\Hnplus[3]\of{\Rn[1]} \sim
            \directsum{\Rnplus[2]}{\Lnplus[4]}$ which has three
            factors, because $\Rnplus[2] \cong
            \directsum{\Lnplus[1]}{\Lnplus[1]}$. If
            $\directsum{\Lnplus[m]}{\Lnplus} \sim
            \directsum{\Hnplus[3]\of{\Rn[1]}}{J}$, it follows that
            $\directsum{\Lnplus[m]}{\Lnplus} \sim
            \directsum{\directsum{\Rnplus[2]}{\Lnplus[4]}}{J}$ as
            well. And since the dimension of $J$ is too small to admit
            any non-Lorentz cones, $\directsum{\Lnplus[m]}{\Lnplus}$
            has a simulacrum consisting of only Lorentz cones. We may
            therefore check these cases by partitioning the dimension
            and looking for sums of Lorentz cones of the appropriate
            Lyapunov rank.
          \end{item}
        \end{itemize}

        When all is said and done, the only counterexample that arises
        is $\directsum{\Rnplus[2]}{\Lnplus[4]} \sim
        \Hnplus[3]\of{\Rn[1]}$ corresponding to $m = 2$ and $n = 4$.
      \end{proof}
    \end{proposition}

    The last case that we will single out is
    $\directsum{\Lnplus}{\Lnplus}$, which typically will have
    simulacra except for a few small values of $n$.

    \begin{lemma}\label{lem:Ln-squared is easy to brute force}
      If $\directsum{\Lnplus}{\Lnplus}$ has symmetric simulacra, then it has
      symmetric simulacra of the form
      \begin{equation*}
        \alpha \Hnplus[3]\of{\Cn[1]}
        \oplus
        \Lnplus[n_{1}]
        \oplus
        \Lnplus[n_{2}]
        \oplus
        \cdots
        \oplus
        \Lnplus[n_{k}]
      \end{equation*}
      where $\alpha \in \set{0,1}$ and $k, n_{1}, n_{2}, \ldots, n_{k}
      \in \Nn[1]$ (as in \cref{thm:two sufficient families}).

      \begin{proof}
        For $n \le 2$, the result holds vacuously because
        $\directsum{\Lnplus[n]}{\Lnplus[n]} \cong \Rnplus[2n]$ has no
        simulacra.

        So suppose $J \sim \directsum{\Lnplus}{\Lnplus}$ with $n \ge
        3$. If $J$ is not of the desired form, then the results in
        \cref{sec:irreducible cones,cor:no multiple complex factors}
        can be used to find another symmetric cone having the same
        signature as $J$ that is of the desired form. We need only
        convince ourselves that this new cone is not isomorphic to
        $\directsum{\Lnplus}{\Lnplus}$. For that, we note that the
        simulacra constructed in \cref{sec:irreducible cones,cor:no
          multiple complex factors} all have more than two Lorentz
        factors. $\directsum{\Lnplus}{\Lnplus}$ cannot be isomorphic
        to a direct sum of three or more factors when $n \ge 3$, so we
        are done.
      \end{proof}
    \end{lemma}

    \begin{proposition}
      Suppose $n \in \Nn[1]$. Then $\directsum{\Lnplus}{\Lnplus}$ has
      symmetric simulacra if and only if $n \notin
      \set{0,1,2,3,5,6,7,11,12,13,18}$.

      \begin{proof}
        Suppose that $n \ge 100$, and use the division algorithm to
        define $m$ and $k$ by $5m + k \coloneqq n$. It follows that $m
        \ge 20$ and $k \le 4$. Using the division algorithm once again
        we may define $r \in \set{0,1,2}$ to be the remainder upon
        dividing $m - k^{2} + 1$ by $3$. Next, define
        \begin{align*}
          \alpha &\coloneqq \frac{m - 4k^{2} + 15k - 14 - r}{3},\\
          \gamma &\coloneqq 2m - 22k - \frac{4m - 16k^{2} - 68 + 5r}{3}.
        \end{align*}
        Having assumed that $n \ge 100$, it is straightforward to see
        that $\alpha, \gamma \ge 0$. Both are integers as
        well: rewriting a bit,
        \begin{equation*}
          \alpha
          =
          \frac{m - k^{2} + 1 - r}{3} + \frac{15k - 3k^{2} - 15}{3},
        \end{equation*}
        which is clearly integral considering how $r$ is defined. And
        likewise,
        \begin{equation*}
          \gamma
          =
          2m - 22k - \frac{4m - 16k^{2} + 4 - 4r}{3} - \frac{72 - 9r}{3}.
        \end{equation*}
        All that remains for the $n \ge 100$ case is to define
        \begin{equation*}
          J
          \coloneqq
          \Lnplus[7m+k]
          \oplus
          \Lnplus[m+3k-4]
          \oplus
          \underbrace{
            \Lnplus[4]
            \oplus
            \cdots
            \oplus
            \Lnplus[4]
          }_{\alpha \text{ times }}
          \oplus
          \underbrace{
            \Lnplus[3] 
            \oplus
            \cdots
            \oplus
            \Lnplus[3]
          }_{r \text{ times}}
          \oplus
          \Rnplus[\gamma].
        \end{equation*}
        This cone shares its signature with
        $\directsum{\Lnplus}{\Lnplus}$, and its first factor is
        $\Lnplus[7m+k] = \Lnplus[n+2m]$, so the two are not
        isomorphic.

        For $n < 100$, we must once again resort to brute-force
        enumeration. From \cref{lem:Ln-squared is easy to brute
          force}, we know however that this can be done (relatively)
        easily using the same strategy we used in \cref{thm:improved
          subproblem reduction}, namely by partitioning the desired
        dimension and computing Lyapunov ranks. One would first
        partition $2n-9$ assuming that $\Hnplus[3]\of{\Cn[1]}$ is
        present in the simulacra, searching for a Lyapunov rank of
        $\lyapunovrank{\directsum{\Lnplus}{\Lnplus}} -
        \lyapunovrank{\Hnplus[3]\of{\Cn[1]}} = n^{2} - n - 15$. Should
        that fail\footnote{Partitioning $2n-9$ is a subproblem of
        partitioning $2n$, if you care to save the result.}, $2n$
        itself can be partitioned in hopes of finding
        $\lyapunovrank{\directsum{\Lnplus}{\Lnplus}}$. In
        \cref{sec:L(n)+L(n) simulacra}, we have provided the simulacra
        for all $n \notin \set{0,1,2,3,5,6,7,11,12,13,18}$. It remains
        only to confirm the absence of simulacra in those few cases.
      \end{proof}
    \end{proposition}
  \end{section}

  \section*{Acknowledgments}{
    The authors thank David Sossa and the University of O'Higgins for
    organizing and sponsoring the 2024 Workshop on Variation Analysis
    and Euclidean Jordan Algebras (VAEJA). The results in
    \cref{sec:reducible cones} were proved as a result of the
    workshop.

    G.B. is a member of the Research Group GNCS (Gruppo Nazionale per
    il Calcolo Scientifico) of INdAM (Istituto Nazionale di Alta
    Matematica). G.B. is supported by the ERC Consolidator Grant
    101085607 through the Project eLinoR.
  }

  \appendix
  \begin{section}{The table of Lyapunov ranks}
    In the absence of any restrictions on $n$, \cref{tab:signatures of
      building blocks} would produce incorrect values for
    $\Lnplus[0]$, $\Hnplus[0]\of{\Cn[1]}$, and
    $\Hnplus[1]\of{\quaternions}$. The theorem of Gowda and Tao (on
    which our table is based) is itself based upon a table at the
    bottom of page 97 in Faraut and Korányi~\cite{faraut-koranyi}. The
    latter states no restrictions on $n$, but that can lead to
    misinterpretation in two ways. First, the rows in the table are
    not necessarily distinct. When $n = 1$ for example, we have
    $\Lnplus = \Hnplus\of{\quaternions}$, and when that happens, the
    choice of the $\Hnplus\of{\quaternions}$ row (as opposed to the
    $\Lnplus$ row) is what produces the wrong Lyapunov rank. Second,
    the table is simply not valid in trivial cases.

    Based on the notation, we find it likely that the Lie algebras in
    Faraut and Korányi's table are derived from Satz~IX.3.3 in Braun
    in Koecher~\cite{braun-koecher}. Braun and Koecher have one case
    for $\Lnplus$ with $n \ge 2$, and then require $n \ge 3$ for
    $\Hnplus\of{\Rn[1]}$, $\Hnplus\of{\Cn[1]}$, and
    $\Hnplus\of{\quaternions}$. As in \cref{thm:direct sum
      decomposition}, this ensures that there is no overlap between
    the families. From this we may draw two conclusions about the
    Faraut and Korányi table:

    \begin{enumerate}
      \begin{item}
        When there is ambiguity, the Lorentz cone row should be used.
      \end{item}
      \begin{item}
        A priori, the table is invalid for $n < 2$.
      \end{item}
    \end{enumerate}

    It is interesting to note however that the entries for
    $\Hnplus[2]\of{\Rn[1]} \cong \Lnplus[3]$, $\Hnplus[2]\of{\Cn[1]}
    \cong \Lnplus[4]$, and $\Hnplus[2]\of{\quaternions} \cong
    \Lnplus[6]$ are consistent. In each case we have a choice of rows
    (and should prefer the row for the Lorentz cone), but up to
    isomorphism, either row gives the same Lie algebras. The entire
    table is therefore free of ambiguity for $n \ge 2$, which is how
    we arrive at the $n \ge 2$ condition on \cref{tab:signatures of
      building blocks}.

    For the remaining cases $n \in \set{0,1}$, the computations are
    trivial and we find it much simpler to use a separate table
    (\cref{tab:rnplus signature}) than it would be to make
    \cref{tab:signatures of building blocks} consistent for all $n \ge
    0$.
  \end{section}

  \begin{section}{Simulacra of $\directsum{\Lnplus}{\Lnplus}$}\label{sec:L(n)+L(n) simulacra}

    \setlength{\tabcolsep}{0.7em}

    \begin{table}[H]
        \footnotesize
        \begin{center}
        \begin{tabular}{rl  rl  rl}
          $n$ & $\operatorname{partition}\of{2n}$
            & $n$ & $\operatorname{partition}\of{2n}$
            & $n$ & $\operatorname{partition}\of{2n}$\\
          \midrule
          4 & $\set{1, 2, 5}$
            & 41 & $\set{2, 4, 23, 53}$
            & 71 & $\set{2, 4, 8, 34, 94}$\\
          8 & $\set{1, 5, 10}$
            & 42 & $\set{1, 2, 3, 3, 19, 56}$
            & 72 & $\set{1, 2, 3, 4, 41, 93}$\\
          9 & $\set{1, 2, 3, 12}$
            & 43 & $\set{2, 5, 23, 56}$
            & 73 & $\set{1, 2, 2, 2, 2, 2, 40, 95}$\\
          10 & $\set{2, 5, 13}$
            & 44 & $\set{1, 37, 50}$
            & 74 & $\set{1, 65, 82}$\\
          14 & $\set{1, 10, 17}$
            & 45 & $\set{1, 3, 30, 56}$
            & 75 & $\set{1, 3, 56, 90}$\\
          15 & $\set{1, 3, 6, 20}$
            & 46 & $\set{5, 29, 58}$
            & 76 & $\set{1, 2, 4, 50, 95}$\\
          16 & $\set{2, 2, 2, 2, 2, 22}$
            & 47 & $\set{1, 2, 2, 2, 26, 61}$
            & 77 & $\set{2, 4, 53, 95}$\\
          17 & $\set{2, 4, 5, 23}$
            & 48 & $\set{1, 2, 2, 2, 6, 18, 65}$
            & 78 & $\set{1, 3, 6, 46, 100}$\\
          19 & $\set{1, 2, 2, 2, 5, 26}$
            & 49 & $\set{2, 2, 4, 26, 64}$
            & 79 & $\set{2, 3, 57, 96}$\\
          20 & $\set{2, 13, 25}$
            & 50 & $\set{1, 2, 10, 20, 67}$
            & 80 & $\set{5, 58, 97}$\\
          21 & $\set{1, 2, 12, 27}$
            & 51 & $\set{2, 3, 33, 64}$
            & 81 & $\set{1, 4, 7, 45, 105}$\\
          22 & $\set{1, 17, 26}$
            & 52 & $\set{2, 41, 61}$
            & 82 & $\set{2, 2, 4, 53, 103}$\\
          23 & $\set{1, 3, 12, 30}$
            & 53 & $\set{2, 5, 31, 68}$
            & 83 & $\set{2, 5, 56, 103}$\\
          24 & $\set{1, 2, 2, 2, 2, 6, 33}$
            & 54 & $\set{1, 3, 4, 30, 70}$
            & 84 & $\set{3, 3, 59, 103}$\\
          25 & $\set{2, 3, 12, 33}$
            & 55 & $\set{2, 2, 2, 5, 26, 73}$
            & 85 & $\set{4, 7, 50, 109}$\\
          26 & $\set{5, 13, 34}$
            & 56 & $\set{10, 29, 73}$
            & 86 & $\set{10, 53, 109}$\\
          27 & $\set{2, 2, 2, 12, 36}$
            & 57 & $\set{2, 2, 2, 36, 72}$
            & 87 & $\set{2, 3, 64, 105}$\\
          28 & $\set{1, 5, 13, 37}$
            & 58 & $\set{1, 50, 65}$
            & 88 & $\set{5, 65, 106}$\\
          29 & $\set{1, 2, 2, 2, 2, 2, 7, 40}$
            & 59 & $\set{1, 3, 42, 72}$
            & 89 & $\set{1, 2, 2, 2, 61, 110}$\\
          30 & $\set{1, 2, 3, 4, 9, 41}$
            & 60 & $\set{1, 2, 2, 2, 2, 33, 78}$
            & 90 & $\set{1, 4, 6, 54, 115}$\\
          31 & $\set{3, 8, 9, 42}$
            & 61 & $\set{4, 5, 34, 79}$
            & 91 & $\set{1, 7, 9, 45, 120}$\\
          32 & $\set{1, 26, 37}$
            & 62 & $\set{1, 2, 3, 4, 33, 81}$
            & 92 & $\set{1, 82, 101}$\\
          33 & $\set{1, 3, 20, 42}$
            & 63 & $\set{1, 2, 48, 75}$
            & 93 & $\set{1, 2, 75, 108}$\\
          34 & $\set{2, 25, 41}$
            & 64 & $\set{1, 2, 2, 44, 79}$
            & 94 & $\set{1, 2, 2, 2, 2, 61, 118}$\\
          35 & $\set{2, 8, 13, 47}$
            & 65 & $\set{4, 7, 34, 85}$
            & 95 & $\set{2, 2, 4, 64, 118}$\\
          36 & $\set{3, 3, 19, 47}$
            & 66 & $\set{1, 2, 4, 5, 33, 87}$
            & 96 & $\set{1, 2, 4, 5, 57, 123}$\\
          37 & $\set{1, 2, 2, 5, 14, 50}$
            & 67 & $\set{3, 3, 4, 37, 87}$
            & 97 & $\set{2, 5, 68, 119}$\\
          38 & $\set{1, 2, 4, 19, 50}$
            & 68 & $\set{1, 2, 7, 38, 88}$
            & 98 & $\set{2, 3, 8, 57, 126}$\\
          39 & $\set{1, 2, 27, 48}$
            & 69 & $\set{1, 2, 2, 6, 37, 90}$
            & 99 & $\set{2, 2, 2, 72, 120}$\\
          40 & $\set{2, 2, 4, 19, 53}$
            & 70 & $\set{3, 3, 47, 87}$
            & 100 & $\set{2, 85, 113}$\\
        \end{tabular}
        \end{center}
        \normalsize
    \end{table}
  \end{section}

  \bibliographystyle{mjo}
  \bibliography{local-references}
\end{document}

%% file: mjo-algebra.tex
\ifx\havemjoalgebra\undefined
\def\havemjoalgebra{1}

\ifx\operatorname\undefined
  \usepackage{amsopn}
\fi

\input{mjo-common} 

\newcommand*{\zero}[1]{ 0_{{#1}} }

\ifdefined\newglossaryentry
  \newglossaryentry{zero}{
    name={\ensuremath{\zero{R}}},
    description={the additive identity element of $R$},
    sort=z
  }
\fi

\newcommand*{\unit}[1]{ 1_{{#1}} }

\ifdefined\newglossaryentry
  \newglossaryentry{unit}{
    name={\ensuremath{\unit{R}}},
    description={the multiplicative identity (unit) element of $R$},
    sort=u
  }
\fi

\newcommand*{\directsum}[2]{ {#1}\oplus{#2} }

\newcommand*{\directsummany}[3]{ \binopmany{\bigoplus}{#1}{#2}{#3} }

\newcommand*{\alg}[1]{\operatorname{alg}\of{{#1}}}
\ifdefined\newglossaryentry
  \newglossaryentry{alg}{
    name={\ensuremath{\alg{X}}},
    description={the (sub)algebra generated by $X$},
    sort=a
  }
\fi

\newcommand*{\ideal}[1]{\operatorname{ideal}\of{{#1}}}
\ifdefined\newglossaryentry
  \newglossaryentry{ideal}{
    name={\ensuremath{\ideal{X}}},
    description={the ideal generated by $X$},
    sort=i
  }
\fi

\newcommand*{\polyring}[2]{{#1}\left[{#2}\right]}
\ifdefined\newglossaryentry
  \newglossaryentry{polyring}{
    name={\ensuremath{\polyring{R}{X}}},
    description={polynomials with coefficients in $R$ and variable $X$},
    sort=p
  }
\fi

\newcommand*{\variety}[1]{ \mathcal{V}\of{{#1}} }
\ifdefined\newglossaryentry
  \newglossaryentry{variety}{
    name={\ensuremath{\variety{I}}},
    description={variety corresponding to the ideal $I$},
    sort=v
  }
\fi

\fi

%% file: mjo-common.tex
\ifx\havemjocommon\undefined
\def\havemjocommon{1}

\ifx\mathbb\undefined
  \usepackage{amsfonts}
\fi

\ifx\restriction\undefined
  \usepackage{amssymb}
\fi

\newcommand*{\of}[1]{ \left({#1}\right) }

\newcommand*{\qty}[1]{ \left({#1}\right) }

\newcommand*{\pair}[2]{ \left({#1},{#2}\right) }

\newcommand*{\lcm}[1]{ \operatorname{lcm}\of{{#1}} }
\ifdefined\newglossaryentry
  \newglossaryentry{lcm}{
    name={\ensuremath{\lcm{X}}},
    description={the least common multiple of the elements of $X$},
    sort=l
  }
\fi

\newcommand*{\restrict}[2]{{#1}{\restriction}_{#2}}
\ifdefined\newglossaryentry
  \newglossaryentry{restriction}{
    name={\ensuremath{\restrict{f}{X}}},
    description={the restriction of $f$ to $X$},
    sort=r
  }
\fi

\newcommand*{\Nn}[1][n]{
  \mathbb{N}\if\detokenize{#1}\detokenize{1}{}\else^{#1}\fi
}

\ifdefined\newglossaryentry
  \newglossaryentry{N}{
    name={\ensuremath{\Nn[1]}},
    description={the set of natural numbers},
    sort=N
  }
\fi

\newcommand*{\Zn}[1][n]{
  \mathbb{Z}\if\detokenize{#1}\detokenize{1}{}\else^{#1}\fi
}

\ifdefined\newglossaryentry
  \newglossaryentry{Z}{
    name={\ensuremath{\Zn[1]}},
    description={the ring of integers},
    sort=Z
  }
\fi

\newcommand*{\Qn}[1][n]{
  \mathbb{Q}\if\detokenize{#1}\detokenize{1}{}\else^{#1}\fi
}

\ifdefined\newglossaryentry
  \newglossaryentry{Q}{
    name={\ensuremath{\Qn[1]}},
    description={the field of rational numbers},
    sort=Q
  }
\fi

\newcommand*{\Rn}[1][n]{
  \mathbb{R}\if\detokenize{#1}\detokenize{1}{}\else^{#1}\fi
}

\ifdefined\newglossaryentry
  \newglossaryentry{R}{
    name={\ensuremath{\Rn[1]}},
    description={the field of real numbers},
    sort=R
  }
\fi

\newcommand*{\Cn}[1][n]{
  \mathbb{C}\if\detokenize{#1}\detokenize{1}{}\else^{#1}\fi
}

\ifdefined\newglossaryentry
  \newglossaryentry{C}{
    name={\ensuremath{\Cn[1]}},
    description={the field of complex numbers},
    sort=C
  }
\fi

\newcommand*{\Fn}[1][n]{
  \mathbb{F}\if\detokenize{#1}\detokenize{1}{}\else^{#1}\fi
}

\ifdefined\newglossaryentry
  \newglossaryentry{F}{
    name={\ensuremath{\Fn[1]}},
    description={a generic field},
    sort=F
  }
\fi

\newcommand*{\binopmany}[4]{
  \mathchoice{ \underset{#2}{\overset{#3}{#1}}{#4} }
             { {#1}_{#2}^{#3}{#4} }
             { {#1}_{#2}^{#3}{#4} }
             { {#1}_{#2}^{#3}{#4} }
}

\fi

%% file: mjo-arrow.tex
\ifx\havemjoarrow\undefined
\def\havemjoarrow{1}

\input{mjo-common} 

\ifx\operatorname\undefined
  \usepackage{amsopn}
\fi

\newcommand*{\const}[1]{\operatorname{const}_{{#1}}}

\ifdefined\newglossaryentry
  \newglossaryentry{const}{
    name={\ensuremath{\const{a}}},
    description={the constant function that always returns $a$},
    sort=c
  }
\fi

\newcommand*{\identity}[1]{ \operatorname{id}_{{#1}} }

\ifdefined\newglossaryentry
  \newglossaryentry{identity}{
    name={\ensuremath{\identity{X}}},
    description={the identity function or arrow on $X$},
    sort=i
  }
\fi

\newcommand*{\inverse}[1]{ #1^{-1} }

\fi

%% file: mjo-cone.tex
\ifx\havemjocone\undefined
\def\havemjocone{1}

\ifx\succcurlyeq\undefined
  \usepackage{amssymb} 
\fi

\input{mjo-common} 
\input{mjo-linear_algebra} 

\newcommand*{\dual}[1]{ #1^{*} }

\newcommand*{\Rnplus}[1][n]{ \Rn[#1]_{+} }

\newcommand*{\Lnplus}[1][n]{ \mathcal{L}^{{#1}}_{+} }

\newcommand*{\Snplus}[1][n]{ \Sn[#1]_{+} }

\newcommand*{\Hnplus}[1][n]{ \Hn[#1]_{+} }

\newcommand*{\posops}[2][]{
  \pi\of{ {#2}
    \if\relax\detokenize{#1}\relax
      {}%
    \else
      {,{#1}}%
    \fi
  }
}

\newcommand*{\lyapunovrank}[1]{ \beta\of{ {#1} } }

\fi

%% file: mjo-linear_algebra.tex
\ifx\havemjolinearalgebra\undefined
\def\havemjolinearalgebra{1}

\ifx\lvert\undefined
  \usepackage{amsmath} 
\fi

\ifx\ocircle\undefined
  \usepackage{wasysym}
\fi

\ifx\clipbox\undefined
  \usepackage{trimclip}
\fi

\input{mjo-common} 

\newcommand*{\ip}[2]{\left\langle{#1},{#2}\right\rangle}

\newcommand*{\spectrum}[1]{\sigma\of{{#1}}}
\ifdefined\newglossaryentry
  \newglossaryentry{spectrum}{
    name={\ensuremath{\spectrum{L}}},
    description={the set of all eigenvalues of $L$},
    sort=s
  }
\fi

\newcommand*{\rref}[1]{\operatorname{rref}\of{#1}}
\ifdefined\newglossaryentry
  \newglossaryentry{rref}{
    name={\ensuremath{\rref{A}}},
    description={the reduced row-echelon form of $A$},
    sort=r
  }
\fi

\newcommand*{\boundedops}[2][]{
  \mathcal{B}\of{ {#2}
    \if\relax\detokenize{#1}\relax
      {}%
    \else
      {,{#1}}%
    \fi
  }
}

\newcommand*{\Sn}[1][n]{ \mathcal{S}^{#1} }
\ifdefined\newglossaryentry
  \newglossaryentry{Sn}{
    name={\ensuremath{\Sn}},
    description={the set of $n$-by-$n$ real symmetric matrices},
    sort=Sn
  }
\fi

\newcommand*{\Hn}[1][n]{ \mathcal{H}^{#1} }
\ifdefined\newglossaryentry
  \newglossaryentry{Hn}{
    name={\ensuremath{\Hn}},
    description={the set of $n$-by-$n$ complex Hermitian matrices},
    sort=Hn
  }
\fi

\newcommand*{\GL}[2]{\operatorname{GL}_{#1}\of{#2}}
\ifdefined\newglossaryentry
  \newglossaryentry{GL}{
    name={\ensuremath{\GL{n}{\mathbb{F}}}},
    description={the general linear group of order $n$ over $\mathbb{F}$},
    sort=GL
  }
\fi

\newcommand*{\Der}[1]{\operatorname{Der}\of{#1}}
\ifdefined\newglossaryentry
  \newglossaryentry{Der}{
    name={\ensuremath{\Der{V}}},
    description={the space of all (linear) derivations on $V$},
    sort=Der
  }
\fi

\newcommand*{\Isom}[2][]{
  \operatorname{Isom}\of{ {#2}
    \if\relax\detokenize{#1}\relax
      {}%
    \else
      {,{#1}}%
    \fi
  }
}
\ifdefined\newglossaryentry
  \newglossaryentry{Isom}{
    name={\ensuremath{\Isom{V}}},
    description={the group of isometries on $V$},
    sort=Isom
  }
  \newglossaryentry{Isom2}{
    name={\ensuremath{\Isom[W]{V}}},
    description={the set of isometries from $V$ to $W$},
    sort=Isom
  }
\fi

\fi

%% file: mjo-font.tex
\ifx\havemjofont\undefined
\def\havemjofont{1}

\usepackage[utf8]{inputenc}

\ifx\mathbb\undefined
  \usepackage{amsfonts}
\fi

\DeclareFontFamily{OT1}{pzc}{}
\DeclareFontShape{OT1}{pzc}{m}{it}{<-> s * [1.15] pzcmi7t}{}
\DeclareMathAlphabet{\mathpzc}{OT1}{pzc}{m}{it}

\fi

%% file: mjo-hurwitz.tex
\ifx\havemjohurwitz\undefined
\def\havemjohurwitz{1}

\newcommand*{\quaternions}{\mathbb{H}}

\ifdefined\newglossaryentry
  \newglossaryentry{quaternions}{
    name={\ensuremath{\quaternions}},
    description={the algebra of quaternions},
    sort=H
  }
\fi

\newcommand*{\octonions}{\mathbb{O}}

\ifdefined\newglossaryentry
  \newglossaryentry{octonions}{
    name={\ensuremath{\octonions}},
    description={the algebra of octonions},
    sort=O
  }
\fi

\fi

%% file: mjo-proof_by_cases.tex
\ifx\havemjoproofbycases\undefined
\def\havemjoproofbycases{1}

\usepackage[loadonly]{enumitem}

\usepackage{calc}

\newlength{\singleblskip}
\newlist{pcases}{enumerate}{1}
\setlist[pcases]{
  label=\textit{Case~\arabic*}:~\protect\thiscase\textit{.},
  ref=\arabic*,
  align=left,
  listparindent=\parindent,
  parsep=\parskip,
  itemsep=\singleblskip}

\fi

%% file: mjo-set.tex
\ifx\havemjoset\undefined
\def\havemjoset{1}

\input{mjo-common} 
\input{mjo-font} 

\ifx\operatorname\undefined
  \usepackage{amsopn}
\fi

\ifx\bigtimes\undefined
  \usepackage{mathtools}
\fi

\newcommand*{\set}[1]{\left\lbrace{#1}\right\rbrace}

\newcommand*{\setc}[2]{\left\lbrace{#1}\ \middle|\ {#2}\right\rbrace}

\newcommand*{\powerset}[1]{\mathpzc{P}\of{{#1}}}
\ifdefined\newglossaryentry
  \newglossaryentry{powerset}{
    name={\ensuremath{\powerset{X}}},
    description={the ``powerset,'' or set of all subsets of $X$},
    sort=p
  }
\fi

\newcommand*{\cartprodmany}[3]{ \binopmany{\bigtimes}{#1}{#2}{#3} }

\fi

%% file: mjo-theorem.tex
\ifx\havemjotheorem\undefined
\def\havemjotheorem{1}

\ifx\theoremstyle\undefined
  \usepackage{amsthm}
\fi

\theoremstyle{plain}
\newtheorem{corollary}{Corollary}

\newtheorem{lemma}{Lemma}
\newtheorem{proposition}{Proposition}
\newtheorem{theorem}{Theorem}

\theoremstyle{definition}
\newtheorem{definition}{Definition}

\theoremstyle{remark}

\fi